\documentclass{article}
\usepackage{amsmath, amssymb, amsthm, amsfonts, amsxtra, latexsym, amscd,
pb-diagram, graphics}
\theoremstyle{plain}

\newcommand{\tx}{\otimes }
\newcommand{\ts}{\oplus}

\newcommand{\Lh}{\frak{L}}

\newcommand{\Fm}{\widetilde{F}}
\newcommand{\Fa}{\breve{F}}

\begin{document}
\renewcommand{\proofname}{$\mathbf{Proof}$}
\renewcommand{\refname}{$\mathbf{References}$}
\pagestyle{headings}
\centerline{$\mathbf{ INTRODUCTION\   TO\ ANN}$-$\mathbf{CATEGORIES}$}

\begin{center}
Nguyen Tien Quang\\
\end{center}
$\mathbf{Abstract.}$ In this paper, we present new concepts of Ann-categories, Ann-functors, and a transmission of the structure of categories based on Ann-equivalences. We build Ann-category $End(\mathcal A)$ of Pic-funtors $F : \mathcal A \longrightarrow \mathcal A' $ and prove that each Ann-category can be faithfully embedded  into an almost strictly Ann-category.\\

\centerline{$\mathbf{1. \ Introduction}$}
Monoidal and symmetric monoidal categories were firstly studied by J. Benabou [1], S.Maclane [5], S.Eilenberg and G.M.Kelly [2] and have more and more applications. Hoang Xuan Sinh [3] has studied Gr-categories and Pic-categories are respectively monoidal and symmetric monoidal categories in which objects are invertable and morphisms are isomorphisms. Beside considering the exact problems, Hoang Xuan Sinh has studied the structures of the two classes of these categories and pointed out the invariants which are specific for each class through the  theory of cohomology of groups.\\
\noindent By another direction of research, M.Laplaza [4] has considered the exactness of distributivity constraints in a category have two symmetric monoidal structures with two operations $\oplus, \otimes.$ In the work of Laplaza, the natural mono-morphisms
\[\lambda = \lambda_{A,B,C} : A\otimes(B\oplus C) \longrightarrow (A\otimes B)\oplus (A\otimes C)\]
\[\rho = \rho_{A,B,C} : (A\oplus B)\otimes C \longrightarrow (A\otimes C)\oplus (B\otimes C)\]
together with the natural isomorphisms of the two symmetric monoidal structures must satisfy 24 commutative diagrams which establish the natural relationship between them.\\
\noindent Combining the two above directions, we consider a class of categories which have distributivity constraints by defining the second operation on Pic-categories to obtain the structures like the structure of rings. The operation together with the distributivity constraints have to satisfy some natural conditions so that we can prove the exact theorem and describe the structure of the class of categories. This paper\footnote{This paper has been published (in Vietnamese) in Vietnam Journal of Mathematics Vol. XV, No 4, 1987} presents the system of axioms of Ann-categories and the main results of the structure transference and the embedding theorem to give a base to the proving the exact theorem and studying the structure of Ann-categories later.\\
\noindent Throughout this paper, for the tensor product of two objects $A$ and $B,$ we write $AB$ instead of $A\otimes B ,$ but for the morphisms we still write $f \otimes g$ to avoid confusion with composition.\\

\centerline{$\mathbf{2.\ The\ system\ of\ axioms\ of\ Ann}$-$\mathbf{categories}$}

\indent To prepare for construction the system of axioms of Ann-categories, we first recall the basic concepts of monoidal categories [see [3], [7]).\\
\indent {$\mathbf{\otimes AU}$-$\mathbf{category.}$} An AU-category (or a monoidal category) is a category $\mathcal A$ together with bifunctor $\otimes : \mathcal A \times \mathcal A \longrightarrow \mathcal A,$ a fixed object $1 \in \mathcal A$ and the natural isomorphisms $a, l, r:$\\
\[a_{X,Y,Z} : X\otimes(Y\otimes Z) \longrightarrow (X\otimes Y)\otimes Z \quad (X, Y, Z \in Ob \mathcal A)\]
\[l_X : 1\otimes X \longrightarrow X ,\quad r_X : X\otimes 1 \longrightarrow X\]
satisfy the following commutative diagrams:
\clearpage
\begin{equation}
{\scriptsize
\divide \dgARROWLENGTH by 2
\begin{diagram}
\node{(X\otimes (Y\otimes Z))\otimes T} \arrow[2]{e,t}{a_{X,Y,Z} \otimes id_T} \node[2]{((X\otimes Y)\otimes Z)\otimes T}\\ 
\node{X\otimes ((Y\otimes Z)\otimes T)} \arrow{n,l}{a_{X,Y\otimes Z, T}}
\node[2]{(X\otimes Y)\otimes (Z\otimes T)} \arrow{n,r}{ a_{X\otimes Y,Z,T}}\\
\node[2]{X\otimes (Y\otimes (Z\otimes T))}\arrow{nw,b}{id_X \otimes a_{Y,Z,T}}\arrow{ne,b}{a_{X,Y,Z\otimes T}}
\end{diagram}}\nonumber\qquad\qquad(2.1)
\end{equation}
\begin{equation}
{\scriptsize
\divide \dgARROWLENGTH by 2
\begin{diagram}
\node{X\otimes (1\otimes Y)} \arrow[2]{e,t}{a_{X,1,Y}} \arrow{se,b}{id_X \otimes l_Y}
\node[2]{(X\otimes 1)\otimes Y} \arrow{sw,b}{r_X \otimes id_Y} \\
\node[2]{X\otimes Y} 
\end{diagram}}\nonumber\qquad\qquad\qquad(2.2)
\end{equation}
$a$ and $(1, l, r)$ are respectively called associativite constraints and unit constraints for the operation $\otimes$.\\
\indent If $\mathcal A$ and $\mathcal A'$ are two $\otimes$AU-categories, an $\otimes$AU-functor $F : \mathcal A \longrightarrow \mathcal A'$ is a functor together with a natural isomorphism $\widetilde F,  \widetilde {F}_{X,Y} : F(X\otimes Y) \longrightarrow FX \otimes FY$ and an isomorphism $F_1 : F1 \longrightarrow 1'$ satisfies the following commutative diagrams:
{\scriptsize
\[\begin{CD}
F(X\otimes (Y\otimes Z)) @ > \widetilde F >> FX\otimes F(Y\otimes Z)@>id \otimes \widetilde F>>FX\otimes (FY\otimes FZ)\\
@VF(a) VV @.     @VV a'V\\
F((X\otimes Y)\otimes Z)@ > \widetilde F >>F(X\otimes Y)\otimes FZ @>  \widetilde {F}\otimes id >> (FX\otimes FY)\otimes FZ
\end{CD}\qquad\qquad(2.3)\]}
{\scriptsize
\[\begin{CD}
F(1\otimes X) @> \widetilde F >> F1\otimes FX\\
@V F(l) VV           @VV  F_1\otimes idV\\
FX @<l'  << 1' \otimes FX\\
\end{CD}\qquad(2.4) \qquad \qquad
\begin{CD}
F(X\otimes 1) @> \widetilde F >> FX\otimes F1\\
@V F(r) VV           @VV id \otimes F_1 V\\
FX @<r'  << FX \otimes 1'\\
\end{CD}\qquad(2.4')\]}

\indent In particular, $F = ( F, \widetilde F)$ is called an $\otimes$-functor. If F satisfies (2.3) [resp. (2.4) and (2.4') ], then we say F is an $\otimes$-A functor [resp. $\otimes$U-functor ].\\
\indent {$\mathbf{\otimes ACU}$-$\mathbf{category.}$} An $\otimes$ACU-category (or a symmetric monoidal category) is an $\otimes$AU-category together with a natural isomorphism
 \[c = c_{X,Y} : X\otimes Y \longrightarrow Y\otimes X\]
 satisfying the condition $c_{X,Y} . c_{Y,X} = id$ and the following commutative diagram:
{\scriptsize
\[\divide \dgARROWLENGTH by 2
 \begin {diagram}
\node[2]{(X\otimes Y)\otimes Z}\arrow{se,t}{c_{X\otimes Y,Z}}\\
\node{X\otimes (Y\otimes Z)}\arrow{ne,t}{a_{X,Y,Z}}\arrow{s,l}{id_{X}\otimes c_{Y,Z}}
\node[2]{Z\otimes (X\otimes Y)}\arrow{s,r}{a_{Z,X,Y}}\\
\node{X\otimes (Z\otimes Y)}\arrow{se,b}{a_{X,Y,Z}}
\node[2]{(Z\otimes X)\otimes Y}\\
\node[2]{(X\otimes Z)\otimes Y}\arrow{ne,b}{c_{X,Z}\otimes id_{Y}}
\end{diagram}\qquad\qquad(2.5)\]}
\indent An $\otimes$ACU-functor $F : \mathcal A \longrightarrow \mathcal A'$ between two $\otimes$-categories is an $\otimes$AU-functor satisfying the following commutative diagrams:
{\scriptsize 
\[
\divide \dgARROWLENGTH by 2
\begin{diagram}
\node{F(X\otimes Y)}\arrow{e,t}{\widetilde F}\arrow{s,l}{F(c)}
\node{FX\otimes FY}\arrow{s,r}{c}\\
\node{F(Y\otimes X)}\arrow{e,t}{\widetilde F}\node{FY\otimes FX}
\end{diagram}\qquad(2.6)\]}
\clearpage
In particular, the pair $(F, \widetilde F)$ satisfying the diagrams (2.3) and (2.6) is called an $\otimes$AC-functor.

\indent {$\mathbf{Proposition\ 2.0.}$} If $F : \mathcal A \longrightarrow \mathcal A' $, $G: \mathcal A' \longrightarrow \mathcal A''$ are $\otimes$AU-functors (resp $\otimes$ACU-functors), the composition $G\circ F$ is also a $\otimes$ AU-functor [resp. $\otimes$ ACU-functor ] together with the natural isomorphism $\widetilde{G\circ F}$ and the isomorphism $\widehat{G\circ F}$ given by the following commutative diagrams:
{\scriptsize \[\begin{CD}
GF(X\otimes Y) @> \widetilde {GF} >> GFX\otimes GFY\\
@|          @AA\widetilde G A\\
G(F(X\otimes Y)) @>  G(\widetilde F) >> G(FX\otimes FY)\\
\end{CD}\qquad (2.7)\qquad \qquad
\begin{CD}
GF1 @> (GF)_1 >> 1''\\
@|          @AA G_1 A\\
G(F1) @>  G(F_1) >> G1'\\
\end{CD}\qquad(2.8)\]}
\indent {$\mathbf{\otimes}$-$\mathbf{morphism.}$} An $\otimes$-morphism between two $\otimes$-functors $(F, \widetilde F) , (G, \widetilde G)$ is the functor morphism $\alpha : F \longrightarrow G$ satisfying the commutative diagram:
{\scriptsize 
\[
\divide \dgARROWLENGTH by 2
\begin{diagram}
\node{F(X\otimes Y)}\arrow{e,t}{\widetilde F}\arrow{s,r}{\alpha_{X\otimes Y}}
\node{FX\otimes FY}\arrow{s,l}{\alpha_X \otimes \alpha_Y}\\
\node{G(X\otimes Y)}\arrow{e,t}{\widetilde G}\node{GX\otimes GY}
\end{diagram}\qquad\qquad(2.9)\]}
\indent{$\mathbf{Pic}$-$\mathbf{category.}$} We say that $\mathcal A$ is a Pic-category if $\mathcal A$ is an $\otimes$ACU-category in which every object is invertable and every morphism is isomorphism.\\
\indent {$\mathbf{Definition\ 2.1.}$} An Ann-category consists of \\
\indent i) The category $\mathcal A$ together with two binary functors $\oplus, \otimes : \mathcal A \times \mathcal A \longrightarrow \mathcal A;$\\
\indent ii) The fixed object $0\in \mathcal A$ together with the natural isomorphism $a^+, c, g, d$ such that ($ \mathcal A, \oplus , a^+, c, (0, g, d)$) is a Pic-category;\\
\indent iii) The fixed object $1 \in \mathcal A$ together with the natural isomorphism $a, l, r$ such that $(\mathcal A, \otimes, a, (1, l, r))$ is a $\otimes$AU-category;\\
\indent iv) The natural isomorphisms $\mathcal L, \mathcal R:$
\[\mathcal L_{A,X,Y} : A \otimes (X\oplus Y) \longrightarrow (A\otimes X) \oplus ( A\otimes Y)\]
\[\mathcal R_{X,Y,A} : (X\oplus Y) \otimes A \longrightarrow (X\otimes A)\oplus (Y\otimes A)\]
such that the following conditions are satisfied:\\
{\it (Ann-1)}  For every object $A\in \mathcal A,$ the pairs $(L^{A}, \breve{L^{A}})$,  $(R^{A}, \breve{R^{A}})$ defined by the equations:
 \[\begin{cases}
L^A : X \longmapsto A\otimes X\\
L^A(u) = id_A \otimes u , u : X\longrightarrow Y\\
{\breve L}^{A}_{X,Y} ={ \mathcal L}_{A,X,Y}
\end{cases}\qquad\qquad
 \begin{cases}
R^A : X \longmapsto X\otimes A\\
R^A(u) = u \otimes id_A  , u : X\longrightarrow Y\\
{\breve R}^{A}_{X,Y} ={ \mathcal R}_{X,Y,A}
\end{cases}\]
are $\otimes$AC-functors.\\
{\it (Ann-2)} For every $A, B, X, Y \in \mathcal A$, the following diagrams are commutative:
{\scriptsize 
\[
\divide \dgARROWLENGTH by 2
\begin{diagram}
\node {A((B X)\oplus (BY))} \arrow {s,l}{\breve{L}^A}
\node{A (B (X\oplus Y))}\arrow{w,t}{id_A\otimes\breve{ L}^B}\arrow{e,t}{a_{A,B,X\oplus Y}} 
 \node{(AB)(X\oplus Y)}\arrow {s,r} {\breve{ L}^{A B}}
\\
\node{(A(B X))\oplus (A(BY))}\arrow[2]{e,t}{a_{A,B,X}\oplus a_{A,B,Y}}\node[2]{((AB)X) \oplus ((A B) Y)}
\end{diagram}\qquad\qquad\qquad(2.10)\]}

{\scriptsize 
\[
\divide \dgARROWLENGTH by 2
\begin{diagram}
\node {(X \oplus Y)(BA)}\arrow{e,t}{a_{X \oplus Y,B,A}}\arrow{s,l}{\breve {R}^{BA}}
\node{((X\oplus Y)B)A}\arrow{e,t}{\breve {R}^{B} \otimes id} 
\node{((XB)\oplus (YB) A}\arrow{s,r}{\breve {R}^A} \\
\node{(X(B A))\oplus (Y(B A))}\arrow[2]{e,t}{a_{X,B,A}\oplus a_{Y,B,A}}
\node[2]{((X B)A)\oplus ((Y B) A)}
\end{diagram}\qquad\qquad\qquad\qquad(2.10')\]}

{\scriptsize 
\[
\divide \dgARROWLENGTH by 2
\begin{diagram}
\node{A((X\oplus Y)B)}\arrow{e,t}{id_A \otimes \breve {R}^{B}}\arrow{s,l}{a_{A,X\oplus Y,B}}
\node{A((X B)\oplus (Y B))}\arrow{e,t}{\breve {L}^{A}}\node{(A(XB))\oplus (A (YB))}\arrow{s,r}{a_{A,X,B}\oplus a_{A,Y,B}}\\
\node{(A(X\oplus Y)) B}\arrow{e,t}{\breve {L}^{A} \otimes id_B}
\node{((AX)\oplus (A Y)) B}\arrow{e,t}{\breve {R}^{B}}\node{((AX) B)\oplus ((AY)B))}
\end{diagram}\qquad\qquad\qquad(2.11)\]}
{\scriptsize 
\[
\divide \dgARROWLENGTH by 2
\begin{diagram}
\node{((A\oplus B)X)\oplus ((A\oplus B)Y)}\arrow{s,l}{\breve{ R}^{X}\oplus \breve {R}^Y}
\node{(A\oplus B)(X\oplus Y)}\arrow{e,t}{\breve {R}^{X\oplus Y}}\arrow{w,t}{\breve {L}^{A\oplus B}}
 \node{(A (X\oplus Y))\oplus (B (X\oplus Y))}\arrow{s,r}{\breve {L}^{A}\oplus \breve {L}^B}\\
\node{((AX)\oplus (B X))\oplus ((AY)\oplus (B X))}\arrow[2]{e,t}{v}\node[2]{((AX)\oplus (A Y)) \oplus ((B X)\oplus (B Y))}
\end{diagram}\qquad(2.12)\]}
in which $v = v_{A,B,C,D} : (A\oplus B)\oplus (C\oplus D) \longrightarrow  (A\oplus C)\oplus (B\oplus D)$ is the functor morphism built uniquely from the morphisms $a^+, c$ and $id$ in Pic-category $(\mathcal A, \oplus)$ and we omit subscripts for convenience.\\
{\it (Ann-3)} For object $1,$ the following diagrams are commutative:
{\scriptsize 
\[
\divide \dgARROWLENGTH by 2
\begin{diagram}
\node{1\otimes (X\oplus Y)} \arrow[2]{e,t}{\breve{L}^{1}_{X,Y}} \arrow{se,b}{l_{X\oplus Y}}
\node[2]{(1\otimes X)\oplus (1\otimes Y)}\arrow {sw,b}{l_X \oplus l_Y} \\
\node[2]{X\oplus Y}
\end{diagram}\qquad\qquad(2.13)\]}
{\scriptsize 
\[
\divide \dgARROWLENGTH by 2
\begin{diagram}
\node{(X\oplus Y)\otimes 1} \arrow[2]{e,t}{\breve{R}^{1}_{X,Y}} \arrow{se,b}{r_{X\oplus Y}}
\node[2]{(X\otimes 1)\oplus (Y\otimes 1)}\arrow {sw,b}{r_X \oplus r_Y} \\
\node[2]{X\oplus Y}
\end{diagram}\qquad\qquad(2.13')\]}

\indent{$\mathbf{Comment\ 2.2.}$} The commutative diagrams (2.10), (2.10'), (2.11) respectively mean that $(a_{A,B,-}) : L^{A}. L^{B} \longrightarrow L^{A\oplus B}$, $(a_{-,A,B}) : R^{A\oplus B} \longrightarrow R^{A}. R^{B}$, $(a_{A,-,B}) : L^{A} . R^{B} \longrightarrow R^{B} .L^{A}$ are $\oplus$ morphisms. The diagram (2.12) says that $(\mathcal L_{-,A,B})$ is $\oplus$-morphism from $R^{A\oplus B}$ to $R^{A}\oplus R^{B}$ and $(\mathcal R_{A,B,-})$ is $\oplus$-morphism from $L^{A\oplus B}$ to $L^{A}\oplus L^{B}$ [the sum of two $\oplus$-functors $(F,\breve F)$ and $(G,\breve G)$ in $\mathcal A$ defined by:
\[(F\oplus G)X = FX\oplus GX \]
\[\breve{F\oplus G} = v\circ(\breve F \oplus \breve G)\qquad (2.14)\]
($v$ is mentioned in definition 2.1)]\\
The diagrams (2.13), (2.13') respectively mean that $l : L^{1} \longrightarrow id_A$ , $r : R^{1} \longrightarrow id_A$ are $\oplus$-morphisms.\\
\indent {$\mathbf{Definition\ 2.3.}$} Let $\mathcal A$ and $\mathcal A'$ be Ann-categories. An Ann-functor $F : \mathcal A \longrightarrow \mathcal A'$ is a functor together with natural isomorphisms $\breve F, \widetilde F$ such that $(F, \breve F)$ is a $\oplus$AC-functor, $(F,\widetilde F)$ is a $\otimes$A-functor. Furthermore, $\breve F$ and $\widetilde F$ are compatible with the distributivity constraints in the sense that the two following diagrams are commutative:
{\scriptsize 
\[
\divide \dgARROWLENGTH by 2
\begin{diagram}
\node{F(X(Y\oplus Z))}\arrow{e,t}{\widetilde {F}}\arrow{s,l}{F(\mathcal {L})}\node{FX F(Y\oplus Z)}\arrow{e,t}{id \otimes \breve {F}}\node{FX(FY\oplus FZ)}\arrow{s,r}{\mathcal {L'}}\\
\node{F((XY)\oplus (XZ))}\arrow{e,t}{\breve{F}}\node{F(X Y) \oplus F(X Z)}\arrow{e,t}{\widetilde{F }\oplus \widetilde {F}}\node{(FX FY)\oplus (FX  FZ)}
\end{diagram}\qquad(2.15)\]}
{\scriptsize 
\[
\divide \dgARROWLENGTH by 2
\begin{diagram}
\node{F((X\oplus Y)\otimes Z)}\arrow{e,t}{\widetilde {F}}\arrow{s,l}{F(\mathcal {R})}
\node{F(X\oplus Y) FZ}\arrow{e,t}{ \breve {F}\otimes id}\node{(FX\oplus FY)FZ}\arrow{s,r}{\mathcal {R'}}\\
\node{F((X Z)\oplus (Y Z))}\arrow{e,t}{\breve{ F}}\node{F(X Z) \oplus F(Y Z)}\arrow{e,t}{\widetilde{F }\oplus \widetilde {F}}\node{(FX FZ)\oplus (FY FZ)}
\end{diagram}\qquad(2.15')\]}
If $F$ is an equivalence, then $F = (F,\breve F,\widetilde F)$ is said to be an Ann-equivalence. In the case $F$ is a faithful functor we say $F = (F,\breve F,\widetilde F)$ is an embedding. We have the following basic result.\\
\indent {$\mathbf{Theorem\ 2.4.}$}\\
\indent {\it i) Each Ann-category can be embedded into an almost strict category in the sense that all of its constraints except for one distributivity constraint (left or right) and commutative constraint, are indentities;\\
\indent ii) The condition $c_{A,A} = id,  A\in \mathcal A$ for the commutative constraints, is neccessary and sufficient to embed $\mathcal A$ into an Ann-category, all of whose constraints except for one distributivity constraint (left or right), are identities}.\\
\indent We left presenting the proof of this theorem to section 5. To do this, we should build quite strict Ann-category $End(\mathcal A)$ (section 5) and structure changing (section 4). The structure changing also plays an important role in later structure considering. The most basic properties of zero-object are presented in section 3.\\
\indent {$\mathbf{Comment.}$} Exactness problems require to prove that in Ann-category $\mathcal A$ every diagram consisting of only morphisms that are built from the isomorphisms $a^+,$ $c,$ $g,$ $d,$ $a,$ $l,$ $r,$ $\mathcal L,$ $\mathcal R$ and the identity depending on the rules $\oplus, \otimes$ are commutative. Thanks to the above embedding theorem, we only need to carry out the proof with quite strict Ann-category. Thus, the work will be much simpler.\\
\indent In fact, we can obtain a better result  than the embedding theorem, that is {\it each Ann-category is also Ann-equivalent to an almost strict Ann-category}. Because proving this is rather long, we choose a presenting through the embedding theorem that leads to the exactness theorem quickly.\\

\centerline{$\mathbf{3.\ Zero}$-$\mathbf{object}$}
\indent Let an Ann-category $\mathcal A= (\mathcal A,\oplus,\otimes)$. For each $A \in \mathcal A$, the pairs $(L^{A}, \Breve {L}^{A})$ and $(R^{A}, \Breve{R}^{A})$ are $\oplus$AC-functors of the Pic-category $(\mathcal A,\oplus)$ so they are $\oplus$ACU-functors (see [3]. I. 4. 8). Thus, we have\\
\indent {$\mathbf{Proposition\ 3.1.}$} {\it In an Ann-category $\mathcal A$, the functors $(L^{A}, \Breve {L}^{A})$ and $(R^{A}, \Breve {R}^{A})$ are $\oplus$ACU-functors. That means there exists an unique isomorphism:
\[ \widehat{L}^{A} : A\otimes 0 \longrightarrow 0 ,\quad  \widehat{R}^{A} : 0\otimes A \longrightarrow 0\]
for every object $A\in \mathcal A$, such that the following diagrams are commutative:}\\
{\scriptsize \[ \begin{CD}
A\otimes X @< L^{A}({g_X}) <<  A\otimes (0\oplus X)\\
@A g_{A\otimes X} AA           @VV \breve L^{A}_{0,X} V \\
0\oplus(A\otimes X) @<\widehat {L}^{A}\otimes id << (A\otimes 0)\oplus (A\otimes X)\\
\end{CD}  \qquad \qquad \qquad
\begin{CD}
A\otimes X@< L^{A}({d_X}) << A\otimes (X\oplus 0)\\
@A d_{A\otimes X} AA           @VV \breve L^{A}_{X,0} V\\
(A\otimes X)\oplus 0 @< id\otimes\widehat L^{A} << (A\otimes X)\oplus (A\otimes 0)\\
\end{CD}\]}
{\scriptsize 
\[\begin{CD}
X \otimes A@< R^{A}({g_X}) <<  (0\oplus X)\otimes A\\
@A g_{X\otimes A} AA           @VV \breve R^{A}_{0,X} V\\
0\oplus(X\otimes A) @<\widehat R^{A}\otimes id <<(0\otimes A)\oplus (X\otimes A)\\
\end{CD}  \qquad \qquad \qquad 
\begin{CD}
X\otimes A@< R^{A}({d_X}) <<  (X\oplus 0)\otimes A\\
@A d_{X\otimes A} AA           @VV \breve R^{A}_{X,0} V\\
(X\otimes A)\oplus 0 @< id\otimes\widehat R^{A} << (X\otimes A)\oplus (0\otimes A)\\
\end{CD}\]}

\indent {$\mathbf{Propositions\ 3.2.}$} {\it In an Ann-category $\mathcal A$, the isomorphisms  $\widehat{L}^{A},  \widehat{R}^{A}$ have the following properties:\\
\indent i) The family $(\widehat{L}-) =\widehat L \ \ (resp.\  the family\  (\widehat{R}-) =\widehat R ) $ is a $\oplus$-morphism from the functor $(R^{0}, \widehat{R}^{0}) \ \ (resp. \ (L^{0}, \widehat{L}^{0}) )$ to the functor $(\theta : A \longmapsto 0, \widehat{ \theta} = {g_{0}}^{-1})$. This means that the following diagrams are commutative:
\[\begin{diagram}
\node{A\otimes 0} \arrow[2]{e,t}{f\otimes id} \arrow{se,b}{\widehat L^{A}}
\node[2]{B\otimes 0}\arrow {sw,b}{\widehat L^{B}} \\
\node[2]{0}
\end{diagram}\qquad \qquad \qquad
\begin{diagram}
\node{0\otimes A} \arrow[2]{e,t}{ id\otimes f} \arrow{se,b}{\widehat R^{A}}
\node[2]{0\otimes B}\arrow {sw,b}{\widehat R^{B}} \\
\node[2]{0}
\end{diagram}\]

{\scriptsize 
\[\begin{CD}
0\oplus 0@> g_0 = d_0 >> 0\\
@A \widehat L^{x}\oplus \widehat L^{y}AA           @AA \widehat L^{x\oplus y}A\\
 (X\otimes 0)\oplus (Y\otimes 0)@< \breve R^{0}<<(X\oplus Y)\otimes 0\\
\end{CD}\qquad \qquad \qquad
\begin{CD}
0\oplus 0@> g_0 = d_0 >> 0\\
@A \widehat R^{x}\oplus \widehat R^{y}AA           @AA \widehat R^{x\oplus y}A\\
 (0\otimes X)\oplus (0\otimes Y)@< \breve L^{0}<<0\otimes (X\oplus Y)\\
\end{CD}\]}\\

\indent ii) For any $A,B \in \mathcal A$, the following diagrams are commutative:
{\scriptsize 
\[
\divide \dgARROWLENGTH by 2
\begin{diagram}
\node{X\otimes (Y\otimes 0)}\arrow{s,l}{a}\arrow{e,t}{id_X \otimes \widehat L^{Y}}
\node{X\otimes 0}\arrow{s,r}{\widehat L^{X}}\\
\node{(X\otimes Y)\otimes 0}\arrow{e,t}{\widehat L^{X\otimes Y}}\node{0}
\end{diagram}\]}
{\scriptsize 
\[
\divide \dgARROWLENGTH by 2
\begin{diagram}
\node{0\otimes (X\otimes Y)}\arrow{s,l}{a}\arrow{e,t}{\widehat R^{X\otimes Y}}\node{0}\\
 \node{(0\otimes X)\otimes Y}\arrow{e,t}{\widehat R^{X}\otimes id_Y}\node{0\otimes Y}\arrow{n,r}{\widehat R^{Y}}
\end{diagram}\]}
{\scriptsize 
\[
\divide \dgARROWLENGTH by 2
\begin{diagram}
\node{X\otimes (0\otimes Y)}\arrow[2]{e,t}{a}\arrow{s,l}{id_X \otimes \widehat R^{Y}}\node[2]{(X\otimes 0)\otimes Y}\arrow{s,r}{\widehat L^{X}\otimes id_Y}\\
\node{X\otimes 0}\arrow{e,t}{\widehat L^{X}}\node{0}\node{0\otimes Y}\arrow{w,t}{\widehat R^{Y}}
\end{diagram}\]}
\indent iii) \[\widehat{L}^{1} = l_0 , \widehat{R}^{1} = r_0\]} .\\
\indent The proofs in detail of the proposition 3.2 and the examples of Ann-category can be seen in [6].\\

\centerline{$\mathbf{4.\ The\ structure\ transference}$}
\indent Let $\mathcal A$ be a categoty and $(\mathcal A',\otimes,a',(1',l',r'))$ be a monoidal category. Let $F : \mathcal A \longrightarrow \mathcal A'$ and $F' : \mathcal A' \longrightarrow \mathcal A$ be the equivalences with isomorphisms $\alpha : F'F \stackrel \sim\longrightarrow id_{\mathcal A} ,\quad \alpha^{'} : FF' \stackrel\sim\longrightarrow id_{\mathcal A'}$. Furthermore, suppose that $\alpha,\alpha^{'}$ satisfy the following conditions:
\[F(\alpha_A) =\alpha^{'}_{FA} ,\quad F'(\alpha^{'}_{A'}) = \alpha_{F'A'}\qquad(4.1)\]
\indent Then, acccording to Hoang Xuan Sinh ([3].I.Section 5), we can define in $\mathcal A$ the operation $\otimes$ thanks to the four-group $(F,F',\alpha,\alpha^{'})$ by the following equations:
\[A\otimes B = F' (FA\otimes FB) \quad A,B \in Ob\mathcal A\]
\[u\otimes v = F' (Fu\otimes Fv)\qquad(4.2)\]
 $u, v$ are morphisms in $\mathcal A.$\\
\indent Furthermore, $F$ and $F'$ become $\otimes$-functors with natural isomorphisms:
\[{\widetilde F}_{A,B} = F (\alpha^{-1}_A \otimes \alpha^{-1}_B) , \quad{\widetilde F'}_{A',B'} = \alpha_{F'A'\otimes F'B'}\qquad(4.3)\]
\indent The pair $(F,\widetilde F)$ induces the unit constraint $(1, l, r)$ in $\mathcal A$ with $1 = F'1'$ and $l, r$ are defined by the following commutative diagrams:
\[\begin{diagram}
\node{1\otimes X}\arrow[2]{e,t}{l_X}\node[2]{X}\\
\node{F'1'\otimes F'FX}\arrow{n,l}{id \otimes \alpha_X}\node{F'(1'\otimes FX)} \arrow{w,t}{\widetilde F'} \arrow{e,t}{F'(l')}\node{F'FX}\arrow{n,r}{\alpha_X}
\end{diagram}\qquad(4.4)\]
\[\begin{diagram}
\node{X\otimes 1}\arrow[2]{e,t}{r_X}\node[2]{X}\\
\node{F'FX\otimes F'1'}\arrow{n,l}{\alpha_X\otimes id}\node{F'(FX\otimes 1')} \arrow{w,t}{\widetilde F'} \arrow{e,t}{F'(r')}\node{F'FX}\arrow{n,r}{\alpha_X}
\end{diagram}\qquad(4.4')\]
\indent The associativity constraint $a$ induced by $(F,\widetilde F)$ is defined by the commutative diagram (2.3). We have\\
\indent {$\mathbf{Proposition\ 4.1.}$} {\it The induced constraints  $a$ and $(1, l, r)$ in $\mathcal A$ are compatible to each other in the sense satisfying the commutative diagram (2.2). Hence, $(\mathcal A,\otimes)$ becomes a monoidal category}.\\
\indent If $\mathcal A'$ is a symmetric monoidal category with the commutative constraint $c'$ then $(\mathcal A,\otimes)$ is also a symmetric monoidal category with the commutative constraint $c$ defined by the commutative diagram (2.6).\\
\indent This leads that if $\mathcal A'$ is an Ann-category then the equivalence $F$ induces a monoidal structure $(\mathcal A,\otimes,a,(1,l,r))$ and a Pic-category structure $(\mathcal A,\oplus,a^+,c,(0,g,d))$. Besides this, the functor $F$ together with the natural isomorphisms $\breve F, \widetilde F$ (defined by the equation (4.3) for each operation $\oplus, \otimes$) permit us to define the distributivity constraints $\mathcal L,\mathcal R$ in $\mathcal A$ by diagrams (2.15), (2.15'). Then, $\mathcal A$ becomes an Ann-category by the following proposition:\\
\indent {$\mathbf{Proposition\ 4.2.}$} {\it Let $\mathcal A$ be a category with two operations $\oplus,\otimes$ in which $(\mathcal A, \oplus)$ is a Pic-category and $(\mathcal A,\otimes)$ is an AU-category. Let $\mathcal A'$ be an Ann-category and $F: \mathcal A \rightarrow \mathcal A'$ is an equivalence together with two natural isomorphisms $\breve{F}$, $\widetilde{F}$ such that $(F,\breve{F})$, $(F,\widetilde F)$ are respectively $\oplus$AC-functor, $\otimes$AU-functor. Then, $\mathcal A$ becomes an Ann-category with natural isomorphisms $\mathcal L, \mathcal R$ defined by the diagrams (2.15), (2.15')}.\\
\indent {$\mathbf{Proof.}$} \\
It is necessary to prove that the isomorphisms $\mathcal L, \mathcal R$ satisfy the commutative diagrams (2.10) - (2.12). For example, we prove that diagram (2.10) is commutative by proving the diagram (4.5) whose parameter is the image of the diagram (2.10) through $F$ is commutative. In this diagram, the parts (I), (VII) are commutative due to the compatibility of $(F, \widetilde F)$ for $a, a'$; the part (II) is commutative by the naturality of $\widetilde F$; the part (III) is commutative by the 

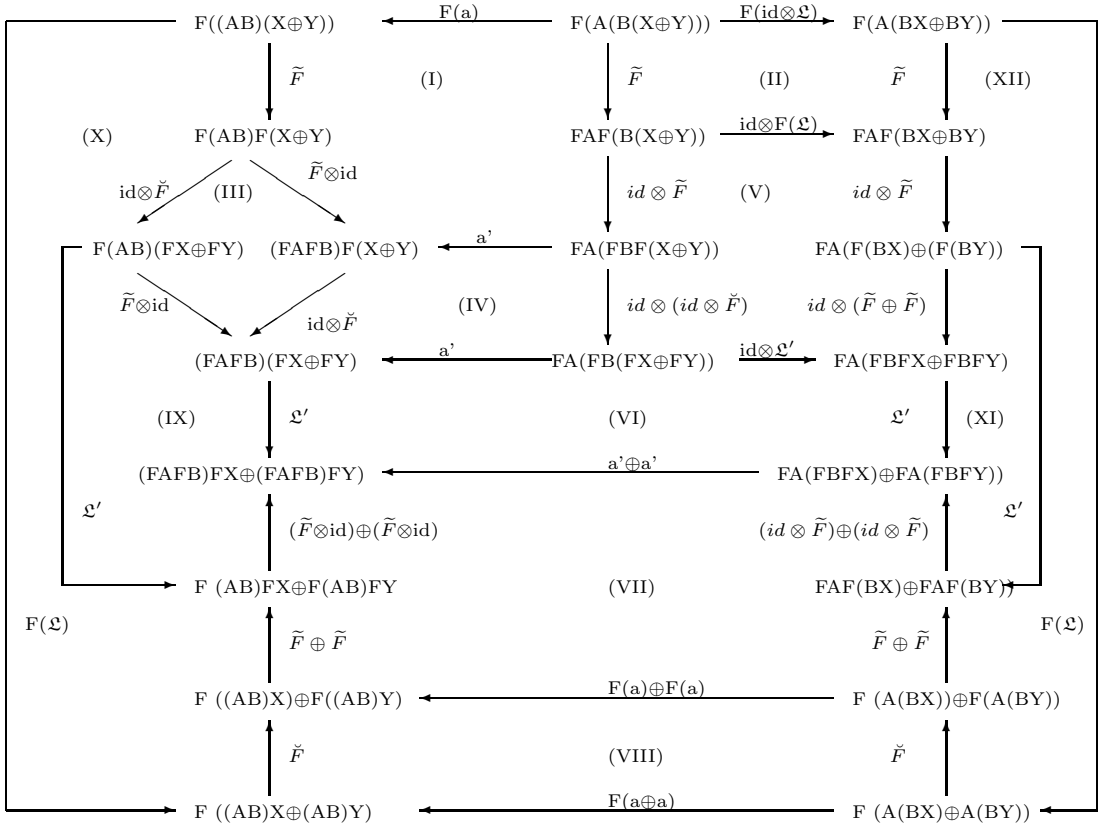
\begin{center}
\setlength{\unitlength}{0.5cm}
\begin{picture}(30,22)
{\scriptsize
\put(5,0){F ((AB)X$\ts$(AB)Y)}                                              \put(22.5,0){F (A(BX)$\ts$A(BY))} 
\put(5,3){F ((AB)X)$\ts$F((AB)Y)}                                           \put(22.5,3){F (A(BX))$\ts$F(A(BY))} 
\put(5,6){F (AB)FX$\ts$F(AB)FY}                                             \put(21.5,6){FAF(BX)$\ts$FAF(BY))} 
\put(3.5,9){(FAFB)FX$\ts$(FAFB)FY)}                                         \put(20.5,9){FA(FBFX)$\ts$FA(FBFY))} 

\put(7,0.6){\vector(0,1){2}}\put(7,3.6){\vector(0,1){2}}\put(7,6.6){\vector(0,1){2}}
\put(25,0.6){\vector(0,1){2}}\put(25,3.6){\vector(0,1){2}}\put(25,6.6){\vector(0,1){2}}

\put(22,0.2){\vector(-1,0){11}}
\put(22,3.2){\vector(-1,0){11}}
\put(20,9.2){\vector(-1,0){10}}

\put(0.5,5){F($\Lh$)}\put(27.5,5){F($\Lh$)}\put(2,8){$\Lh'$}\put(26.5,8){$\Lh'$}

\put(16,1.5){(VIII)}\put(16,6){(VII)}
\put(11,19.5){(I)}\put(20,19.5){(II)}
\put(2,18){(X)}\put(26,19.5){(XII)}
\put(4,10.5){(IX)}\put(16,10.5){(VI)}\put(25.5,10.5){(XI)}
\put(5.5,16.5){(III)}\put(19.5,16.5){(V)}
\put(12,13.5){(IV)}

 \put(16,0.3){F(a$\ts$a)}   \put(16,3.3){F(a)$\ts$F(a)}       \put(16,9.3){a'$\ts$a'}
\put(5,12){(FAFB)(FX$\ts$FY)} \put(14.5,12){FA(FB(FX$\ts$FY))}  \put(22,12){FA(FBFX$\ts$FBFY)} 
\put(2.3,15){F(AB)(FX$\ts$FY)}\put(7,15){(FAFB)F(X$\ts$Y)} \put(15,15){FA(FBF(X$\ts$Y))}  \put(21.5,15){FA(F(BX)$\ts$(F(BY))}
\put(5,18){F(AB)F(X$\ts$Y)} \put(15,18){FAF(B(X$\ts$Y))}  \put(22.5,18){FAF(BX$\ts$BY)} 
\put(5,21){F((AB)(X$\ts$Y))} \put(15,21){F(A(B(X$\ts$Y)))}  \put(22.5,21){F(A(BX$\ts$BY))} 

\put(7,11.6){\vector(0,-1){2}}\put(7,20.6){\vector(0,-1){2}}

\put(25,11.6){\vector(0,-1){2}}\put(25,14.6){\vector(0,-1){2}}\put(25,17.6){\vector(0,-1){2}}\put(25,20.6){\vector(0,-1){2}}
\put(16,14.6){\vector(0,-1){2}}\put(16,17.6){\vector(0,-1){2}}\put(16,20.6){\vector(0,-1){2}}

\put(14.5,12.2){\vector(-1,0){4.5}}\put(14.5,15.2){\vector(-1,0){3}}\put(14.5,21.2){\vector(-1,0){4.5}}
\put(19.5,12.2){\vector(1,0){2}}\put(19,18.2){\vector(1,0){3}}\put(19,21.2){\vector(1,0){3}}

\put(6,17.5){\vector(-3,-2){2.5}}\put(6.5,17.5){\vector(3,-2){2.5}}
\put(3.5,14.5){\vector(3,-2){2.5}}\put(9,14.5){\vector(-3,-2){2.5}}

\put(2,15.2){\line(-1,0){0.5}}\put(1.5,15.2){\line(0,-1){9}}\put(1.5,6.2){\vector(1,0){3}}
\put(27,15.2){\line(1,0){0.5}}\put(27.5,15.2){\line(0,-1){9}}\put(27.5,6.2){\vector(-1,0){1}}

\put(4.5,21.2){\line(-1,0){4.5}}\put(0,21.2){\line(0,-1){21}}\put(0,0.2){\vector(1,0){4.5}}
\put(26.5,21.2){\line(1,0){2.5}}\put(29,21.2){\line(0,-1){21}}\put(29,0.2){\vector(-1,0){1.5}}

  \put(11.5,12.3){a'} \put(12.5,15.3){a'} \put(11.5,21.3){F(a)}
\put(19.5,12.3){id$\tx$$\Lh'$} \put(19.5,18.3){id$\tx$F($\Lh$)} \put(19.5,21.3){F(id$\tx$$\Lh$)}

\put(3,16.5){id$\tx$$\Fa$}\put(8,17){$\Fm$$\tx$id}
\put(3,13.5){$\Fm$$\tx$id}\put(8,13){id$\tx$$\Fa$}

\put(7.5,1.5){$\Fa$}\put(23.5,1.5){$\Fa$}
\put(7.5,4.5){$\Fm\ts\Fm $}\put(23,4.5){$\Fm\ts\Fm $}
\put(7.5,7.5){($\Fm$$\tx$id)$\ts$($\Fm$$\tx$id)}\put(20,7.5){($id\tx\Fm$)$\ts$($id\tx\Fm$)}
\put(7.5,10.5){$\Lh'$}\put(23.5,10.5){$\Lh'$}
\put(7.5,19.5){$\Fm$}\put(16.5,19.5){$\Fm$}\put(23.5,19.5){$\Fm$}

\put(16.5,16.5){$id\tx\Fm$}\put(22.5,16.5){$id\tx\Fm$}
\put(16.5,13.5){$id\tx(id\tx\Fa)$}\put(21.3,13.5){$id\tx(\Fm\ts\Fm)$}

}
\end{picture}
\end{center}
composition of the morphisms; the part (IV) is commutative by the naturality of $a'$; the parts (V), (X), (XII) are commutative  by the definition of $\mathcal L$; the part (VI) is commutative since $\mathcal L'$ satisfies diagram (2.10); the part (VIII) is commutative by the naturality of $\mathcal L'$. Thus, it is true that the parameter is commutative. Because  $F$ is an equivalence we obtain the commutativity of the diagram (2.10).\\
\indent The whole structure transference presented above can be summarized in the following theorem:\\
\indent {$\mathbf{Theorem\ 4.3.}$} {\it Let $\mathcal A$ be a category  and $\mathcal A'$ be an Ann-category . Let $F : \mathcal A \longrightarrow \mathcal A'$ be a category equivalence. Then if $(F,F',\alpha,\alpha')$ is a quadruple satisfying the condition (4.1) then $\mathcal A$ becomes an Ann-category with the induced operations $\oplus, \otimes$ defined by the equation (4.2) and the induced constraints defined by the commutative diagrams (2.3), (2.6), (4.4), (4.4'), (2.15), (2.15') in which $0 = F'0' ; 1 = F'1'$. Furthermore, $F$ is an Ann-equivalence with the natural isomorphisms $\breve{F}$, $\widetilde {F} $ defined by the equation (4.3)}.\\

\indent From the proposition 4.2 we obtain the first result in `stricticizing' the constraints.\\
\indent {$\mathbf{Proposition\ 4.4.}$} {\it Any Ann-category is Ann-equivalent to an Ann-category having AU-strict constraints to the operation $\oplus$}.\\
\indent {$\mathbf{Proof.}$} \\
Assume $\mathcal A$ be an Ann-category. Then, $(\mathcal A,\oplus)$ is ACU-equivalent to a Pic-category $(\mathcal A',\oplus)$ which has AU-strict constraints (see [7], I.2.2.1). By the structure transference we can equip an operation $\otimes$ with the accompanying constraints so $(\mathcal A',\oplus,\otimes)$ becomes an Ann-category that is equivalent to $\mathcal A$.\\

\centerline{$\mathbf{5.\ The\ Ann}$-$\mathbf{category}$ $End(\mathcal A)$ $\mathbf{and\ the\ proof\ of\ the\ embedding\ theorem}$}
\indent Because of the result of the proposition 4.4, in this section we can suppose that $\mathcal A$ is a Pic-category with the AU-strict constraints. Consider the category $End(\mathcal A)$ with the objects are AC-functors $(F,\Breve{ F})$ of $\mathcal A$ ($F$ is an ACU-functor) and the morphisms are the $\oplus$-morphisms. If $F = (F,\Breve{ F}) , G = (G,\Breve{G})$ are two AC-functors then the sum $F \oplus G$ (defined by (2.14)) is also an AC-functor. Hence, $End (\mathcal A)$ becomes a $\oplus$-category when we define the sum of the morphisms $\varphi, \psi$ in $End(\underline A)$ as follows
\[(\varphi + \psi)X = \varphi X \oplus \psi X , \quad X \in \mathcal A.\]
\indent $End (\mathcal A)$ is also a $\otimes$-category with the operation naturally defined by the composition of functors. By these operations we obtain the following theorem which can be proved directly but quite complicated.\\
\indent {$\mathbf{Theorem\ 5.1.}$} {\it $End(\mathcal A)$ is an almost strict Ann-category (in the sense of theorem 2.4) in which\\
\indent i) The zero object is the functor $(\theta : X\longmapsto 0 ,  \Breve{\theta} = id);$\\
\indent ii) The unit object 1 is the identity functor $(I = id_A ,  \Breve{I} = id);$\\
\indent iii) The commutative constrant 
\[c^{\ast} = c^{\ast}_{F,G} : F \oplus G \longrightarrow G\oplus F\] 
defined as follows
\[(c^{\ast}_{F,G}X = c_{FX,GX},\quad X\in \mathcal A)\]
\indent iv) The left distributivity constraint 
\[{\mathcal L}^{\ast} = {\mathcal L}^{\ast}_{F,G,H} : F(G\oplus H) \longrightarrow FG \oplus FH\]
defined by  
\[ ({\mathcal L}^{\ast}_{F,G,H})X = F_{GX,HX} , \quad X\in Ob(\mathcal A)\].}
\indent {$\mathbf{Proposition\ 5.2.}$} {\it Suppose that $\mathcal A$ is an Ann-category with AU-strict constraints for the operation $\oplus$. Then, the correspondence 
\[\begin{matrix} 
\Lambda : &\mathcal A &\longrightarrow& End(\mathcal A) \\
&A&\longmapsto& (L^{A},{\overline L}^{A})\\
&f&\longmapsto& f^{\ast} = (f\otimes id_A) , \quad A\in \mathcal A
\end{matrix}\]
is a faithful functor and an Ann-functor with the natural isomorphisms:
\[\Breve{\Lambda}_{A,B} = ({\mathcal R}_{A,B,X}) : L^{A\oplus B} \longrightarrow {L^{A}} \oplus {L^{B}}\]
\[\widetilde \Lambda_{A,B} = ({a^{-1}}_{A,B,X}) : L^{A\otimes B} \longrightarrow {L^{A}} \otimes {L^{B}} , \quad X\in Ob\mathcal A\]
It means that $\Lambda$ is an embedding $\mathcal A$ into $End(\mathcal A)$}\\
{$\mathbf{Proof.}$} Suppose that $\Lambda(f) =\Lambda(g)$. It follows that $f\otimes id_1 = g\otimes id_1$. Because $1$ is a regular object for the operation $\otimes$, the correspondence $h\longmapsto h\otimes id_1$ is a bijection (object $A$ is regular if the functors $X\longmapsto X\otimes A$ and $X\longmapsto A\otimes X$ are category equivalence). Thus, $f =g$. This means $\Lambda$ is a faithful functor.\\
\indent According to the comment 2.2, ${\widetilde \Lambda}_{A,B} = ({a^{-1}}_{A,B,-})$ and ${\Breve\Lambda}_{A,B} = ({\mathcal R}_{A,B,-})$ are $ \oplus$-functor morphisms.\\
\indent The readers can verify themself that $(\Lambda,\breve \Lambda,\widetilde \Lambda)$ is an Ann-functor directly.\\
\clearpage
\indent {$\mathbf{The\ proof\ of\ the\ theorem\ 2.4}$}\\
\indent i) Suppose that $\mathcal A$ is an Ann-category. From the proposition 4.4, $\mathcal A$ is equivalent to the Ann-category $\mathcal A'$ which has the AU-strict constraints for the addition $\oplus$ through an Ann-equivalence $F.$ Applying the proposition 5.2, we can embed $\mathcal A$ into $End(\mathcal A').$ Hence, the composition
 \[\begin{diagram} \node{\mathcal A}\arrow{e,t}{F}\node{\mathcal A'}\arrow{e,t}{\Lambda} \node{End(\mathcal A')}\end{diagram}\] is the embedding $\mathcal A$ into $End(\mathcal A')$.\\
\indent ii) The sufficient condition: If the Ann-category $\mathcal A$ has a commutative constraints $c$ satisfying the condition $c_{X,X} = id$ then the above Ann-category $\mathcal A'$ will have the commutative constraint $c'=id$ (see [8]). It follows that the commutative constraint of the category $End(\mathcal A')$ (in the proposition 5.1) is also an identity. It is the thing that needs proving.\\
\indent The necessary condition is deduced by the following commutative diagram:\\
\[\begin{CD}
{\mathcal F}(X\oplus X) @> \overline F >> {\mathcal F}(X) \oplus {\mathcal F}(X)\\
@V {\mathcal F}(c_{X,X})VV          @VV c' =id V\\
{\mathcal F}(X\oplus X) @>  \overline F  >> {\mathcal F}(X) \oplus {\mathcal F}(X)\\
\end{CD}\]
in which $(\mathcal F,\overline {\mathcal F},\widetilde {\mathcal F})$ is a faithful Ann-functor which embeds $\mathcal A$ into $\mathcal A'$. Because $\mathcal F$ is faithful so from $\mathcal F(c_{X,X}) = id$, we have $c_{X,X} = id$. The theorem 2.4 has been completely proved.\\
\indent{$\mathbf{Comment.}$} The `strictizing' the constraints in the theorem 2.4 is `maximum'. That means in the general form, we can not embed an Ann-category $\mathcal A$ into an Ann-category $\mathcal A'$ whose constraints are strict except for the commutative constraint. Indeed, if $F: \mathcal A \longrightarrow \mathcal A'$ is the embedding and all of the constraints of $\mathcal A'$ except for the commutative constraint are strict then from the diagram (2.12) we obtain $\nu' = id$. By the definition of the isomorphism $\nu'$ we obtain the isomorphism $c'=id$. Then, applying the theorem 2.4 (ii) we have $c_{X,X} = id$. But that is not true generally.\\
\indent I would like to express my gratefulness to professor Hoang Xuan Sinh for her kindness as giving me guide ideas in writing this paper.

\newpage
\begin{center}

\end{center}
Hanoi University of Education\\
E-mail adress: nguyenquang272002@gmail.com
\end{document}